 
\def\C{{\cal C}}
\def\L{{\cal L}}
\def\M{{\cal M}}
\baselineskip=14pt
\parskip=10pt
\def\Tilde{\char126\relax}

\font\eightrm=cmr8  
\font\eighttt=cmtt8
\magnification=\magstephalf

\parindent=0pt
\overfullrule=0in
\bf
\centerline{The Goulden-Jackson Cluster Method: Extensions, Applications
and Implementations}

\medskip
\it
\centerline{John NOONAN$^1$ and Doron ZEILBERGER\footnote{$^1$}
{\baselineskip=9pt
\eightrm  \raggedright
Department of Mathematics, Temple University,
Philadelphia, PA 19122.
{\baselineskip=9pt
\eighttt [noonan,zeilberg]@math.temple.edu ; 
\break http://www.math.temple.edu/\Tilde 
[noonan,zeilberg]. }
The work of the second author was
supported in part by the NSF. May 2, 1997.
}}
\medskip
\rm
{\bf Abstract:} The powerful 
(and so far under-utilized) Goulden-Jackson Cluster method for finding
the generating function for the number of words avoiding, as factors,
the members of
a prescribed set of `dirty words', is tutorialized and extended in
various directions. 
The authors' Maple implementations, contained in several
Maple packages available from this paper's website
{\tt http://www.math.temple.edu/\Tilde zeilberg/gj.html},
are described and explained.
 
{\bf Preface}
 
In New York City there is a hotel called ESSEX. Once in a while
the bulbs of the first two letters of its neon sign go out,
resulting in the wrong message. This motivates the following
problem. Given a finite alphabet, and a finite set (lexicon) of
`bad words', 
find the number of 
$n$-lettered words in the alphabet that avoid as {\it factors}
(i.e. strings of consecutive letters) any of the dirty words.
More generally, count the number of such words with a prescribed
number of occurrences of obscenities (the previous case being
$0$ bad words), and even more generally, count how many
words are there with a prescribed number of occurrences of
each letter of the alphabet, and a prescribed number of occurrences
of each of the bad words.
 
Many problems in combinatorics, probability, statistics, computer science,
engineering, and  the natural and social sciences, are special cases of, or
can be formulated in terms of, the above scenario.
It is a rather well-kept secret that there exists a powerful
method, the Goulden-Jackson Cluster method[GoJ1][GoJ2], to tackle it.
 
In this paper we start with a motivated and
accessible account
of the method, and then we generalize it in various directions.
Most importantly, we describe our Maple
implementations of both the
original method and of our various extensions. These packages
are obtainable, free of charge, from this paper's very own
website
{\tt http://www.math.temple.edu/\Tilde zeilberg/gj.html }.
 
The Goulden-Jackson Cluster method is very similar, and in some
sense, a generalization of, the method of Guibas and Odlyzko[GuiO],
whose main motivation was Penney-ante games. 
However, philosophically, psychologically, and conceptually, the 
Goulden-Jackson and Guibas-Odlyzko
methods are quite distinct, and we find that the former
is more suitable for our purposes.
 
{\bf The Naive Approach}
 
Before describing the Cluster method, let's review the naive approach.
First, some notation.
Given a word $w=w_1 , \dots , w_n$, a {\it factor} (burrowing the
term from the theory of formal languages) is any of the ${{n+1} \choose {2}}$
words $w_i w_{i+1} \dots  w_{j-1} w_j$, for $1 \leq i \leq j \leq n$.
For example the factors of {\it JOHN} are {\it J, O, H, N, JO, OH, HN,
JOH,OHN, JOHN}, while the factors of {\it DORON} are
{\it D, O, R, O, N, DO, OR, RO, ON, DOR, ORO, RON, DORO, ORON, DORON}.
Note that a given word may occur several times as a factor, for example
the one-letter word {\it O} in {\it DORON}, or the two-letter
words {\it CA} and {\it TI}
in {\it TITICACA}. Also as in formal languages, given
an alphabet $V$, we will denote
the set of all possible words in $V$ by $V^*$.
 
Consider a finite alphabet $V$ with $d$ letters, and 
suppose that we want to keep
track of all factors of length $\leq R+1$, including individual
letters. For every word $w$, of length $\leq R+1$, introduce
a variable $x[w]$. All the $x[w]$ commute with each other.

Define a weight on words $w=w_1, \dots , w_n$, by:
$$
Weight(w)=\prod_{r=1}^{R+1} \prod_{i=1}^{n-r+1} x[w_i, \dots , w_{i+r-1}] 
\quad 
.
$$
For example, if $R=2$, then
$Weight(SEXY)=x[S]x[E]x[X]x[Y]x[SE]x[EX]x[XY]x[SEX]x[EXY]$.
The weight of a set of words (`language') $\L$, $Weight(\L)$, is
defined as the sum of the weights of all the words
belonging to that language. Also, given a language
$\L$ and a letter $v$, we will denote by $\L v$ the set of words obtained
from $\L$ by appending $v$ at the end of each of the words of
$\L$. Thus if $\L=\{SEX, LOON\}$, and $R=1$, then
$Weight(\L)=x[S]x[E]x[X]x[SE]x[EX]+x[L]x[O]^2x[N]x[LO]x[OO]x[ON]$,
and $\L Y=\{SEXY, LOONY\}$.
 
The {\it generating function} 
$$
\Phi_R:=\sum_{w \in V^*} Weight(w) \quad,
$$
stores all the information about the number of words with a prescribed
number of factors of length $\leq R+1$. So, the number of words in
$V^{*}$ that have exactly $n_u$ factors that are $u$ 
for each $u \in V^*$ of $length(u)\leq R+1$,  is the coefficient 
in $\Phi_R$ of the monomial $\prod x[u]^{n_u}$,
where the product extends over the set $\{u \in V^*, length(u)\leq R+1 \} $.
 
If we want the generating function for the number of words
with a prescribed number of bad words and a prescribed number of letters,
we may first compute $\Phi_R$, (where $R+1$ is the maximum length of
a bad word),
and then set $x[v]=s$ for each letter
$v \in V$,  and $x[w]=t$,  if $w$ is a bad word, and $x[w]=1$ otherwise.
The coefficient of $s^n t^m$ in the resulting generating function would
be the number of $n$-letter words with exactly $m$ instances of
bad words occurring as factors.
If we want the generating function for words with no occurrences of
dirty words as factors, we set $t=0$.
 
How to compute $\Phi_R$? For each word $v \in V^*$, of length
$R$, let $Sof[v]$ be the subset of $V^*$ of words that ends with $v$.
Write $v=v_1, \dots, v_R$.
Every word in  $Sof[v]$ is either $v$ itself or of length $> R$,
in which case
chopping the last letter results in an element of $Sof[u]$, for one
of the $d$ $u$'s of the form $i, v_1, \dots , v_{R-1}$.
In symbols
$$
Sof[v]=\{v\} \bigcup_{i \in V} Sof[i,v_1, \dots , v_{R-1}] v_R \quad.
\eqno(SetEq)
$$
 
Since, for any word $w=w_1, \dots, w_n \in V^*$, of length $>R$,
$$
Weight(w_1, \dots , w_n)=Weight(w_1, \dots , w_{n-1}) \cdot  
\prod_{r=1}^{R+1}
x[w_{n-r+1}, \dots, w_n] \quad ,
$$
the system of set equations $(SetEq)$ translates to the linear system
of (algebraic) equations
$$
Weight(Sof[v])=
Weight(v)+ \left ( \prod_{r=1}^{R}
x[v_{R-r+1}, \dots, v_R] \right )
\sum_{i \in V} 
x[i,v_1, \dots, v_R] 
Weight(Sof[i,v_1, \dots , v_{R-1}]) \quad.
\eqno(Linear\_Algebra\_Eq)
$$
 
We have a system of $d^R$ linear equations for $d^R$ unknowns
$Weight(Sof[w])$, $w \in V^{*}, length(w)=R$,
that obviously has a unique solution (on combinatorial grounds!).
Since the coefficients are polynomials 
(in fact monomials) in the variables
$x[w]$, $w \in V^{*}, length(w) \leq R+1$, the solutions
$Weight(Sof[v])$ must be rational functions in these variables.
 
After solving the system, we get $\Phi_R$ from
$$
\Phi_R= \sum_{w \in V^* \,, \, length(w)<R} Weight(w)
+\sum_{w \in V^{*}\, , \, length(w)=R} Weight(Sof[w]) \quad.
$$
Since the first sum is a polynomial and
the second sum is a finite sum of rational functions, it follows that
$\Phi_R$ is a rational function. 
Hence every specialization, as described above, is
also a rational function of its variables.
 
The Maple Implementation of the Naive Approach is contained in
the package {\tt NAIVE}. After downloading it from
this paper's webpage to your working directory,
go into Maple by typing {\tt maple}, followed by {\tt [Enter]}.
Once in Maple, load the package by typing {\tt read NAIVE;}.
To get on-line help, type {\tt ezra();}, for a list of the
procedures, and {\tt ezra(procedure\_name);}, for instructions
how to use a specific function. The most important function
is {\tt PhiR} that computes $\Phi_R$. The function call is
{\tt PhiR(Alphabet, R,x)}, where {\tt Alphabet} is the {\it set}
of letters,
{\tt R} is the {\it non-negative integer} $R$, and
{\tt x} is the {\it variable-name} for the indexed
variables $x[w]$. For example, {\tt PhiR($\{1\}$,0,z);} should give
{\tt 1/(1-z[1]) }.
 
The other procedures are {\tt Naivegf}, {\tt Naivest}, and
{\tt Naives}, that compute, {\it the long way}, what the
procedures {\tt GJgf}, {\tt GJst} and {\tt GJs}
of the package {\tt DAVID\_IAN}, to be described shortly,
compute fast. Their
main purpose is to check the validity of {\tt DAVID\_IAN} and the
other packages described later in this paper. The readers are warned
only to use them for curiosity.

{\bf The Drawback of the Naive Approach}
 
In order to get the generating function 
$\sum_{n=0}^{\infty} a(n)s^n$, where
$a(n):=$ number of words in $\{A, \dots , Z \}^*$ of length
$n$ with no SEX in it (as a factor), we need to solve  a system of
$26^2$ equations and $26^2$ unknowns, then plug in
$x[A]= \dots =x[Z]=s$, $x[AA]= \dots =x[ZZ]=1$,
$x[AAA]= \dots =x[ZZZ]=1$, except for $x[SEX]=0$.
For some economy, we could have made the substitution at the
equations themselves, before solving them,
but we would still  have to solve a system of that size.
 
If we wanted to find the generating function for SEX-less words
with an arbitrary size alphabet, then the above method is
not even valid in principle. Luckily, we have the
powerful {\it Goulden-Jackson Cluster method}, that can handle
such problems very efficiently.
 
{\bf The Most Basic Version of the Goulden-Jackson Cluster Method}
 
Consider a finite alphabet $V$, and a finite set of {\it bad words},
$B$. It is required to find $a(n):=$ the number of words of length
$n$ that {\it do not contain}, as factors, any of the members of the
set of bad words $B$. 
For example if $V=\{ E, S, X \}$, and
$B=\{ SEX, XE \}$, then $a(0)=1, a(1)=3, a(2)=8, a(3)=20, \dots$.
 
Of course we may assume that any factor 
of a bad word of $B$ is not in $B$, since then the longer word
would be superfluous, and can be deleted from the set of banned words.
For example it is not necessary to ban both {\it SEX} and {\it SEXY},
since any word that contains {\it SEXY} in it would also contain
{\it SEX}, and hence the set of words avoiding {\it SEX} and
{\it SEXY} is identical to the set of words avoiding
{\it SEX}.
 
As is often the case in combinatorics, we compute the generating
function $f(s)= \sum_{n=0}^{\infty} a(n)s^n$ rather than $a(n)$
directly. We know from the Naive section that this is a rational
function of $s$, but this fact will emerge again from the Cluster
algorithm, and this time the algorithm is efficient.
 
The methodology is the venerable Inclusion-Exclusion paradigm
that, depending on one's specialty, is sometimes known as
{\it M\"obius Inversion} and {\it Sieve methods}. The essence
of the method is to replace the straight counting of
a hard-to-count set of `good guys' by the weighted count
of the much larger set of pairs 
 
{\it [arbitrary guy, arbitrary subset of his sins]}, 
 
where the 
weight is $(-1)$ to the power of the cardinality of the subset
of his sins.
 
Introducing the weight on words $weight(w):=s^{length(w)}$,
$f(s)$, is the weight enumerator of the set of words, $\L(B)$ that
avoid the members of $B$ as factors, i.e.
$$
f(s)= \sum_{w \in \L(B)} weight(w) \quad .
$$
The trick is to add $0$ to both sides and rewrite this as
$$
f(s)= \sum_{w \in V^*} weight(w) 0^{\,[\,number \,\, of \,\,factors \,\, 
of \,\,w \,\, that \,\,belong \,\,to\,\, B\,]}
 \quad ,
$$
and then use the following deep facts:
$$
0= 1+ (-1) \quad ,
\eqno(i)
$$
$$
0^r= \cases{1,& if $r=0$;\cr
            0,& if $r>0$.\cr} 
\eqno(ii)
$$
and for any finite set $A$,
$$
\prod_{a \in A} 0 = \prod_{a \in A} (1+(-1))=
\sum_{S \subset A} (-1)^{|S|} , \quad
$$
where as usual, $|S|$ denotes the cardinality of $S$.
 
We now have,
$$
f(s)= \sum_{w \in V^*} weight(w) 0^{\,[\,number \,\, of \,\,factors \,\, 
of \,\,w \,\, that \,\,belong \,\,to\,\, B\,]}=
$$
$$
\sum_{w \in V^*} weight(w) (1+(-1))
^{\,[\,number \,\, of \,\,factors \,\, 
of \,\,w \,\, that \,\,belong \,\,to\,\, B\,]}=
$$
$$
\sum_{w \in V^*} 
\sum_{S \subset Bad(w)} (-1)^{|S|} s^{length(w)} \quad ,
$$
where $Bad(w)$ is the set of factors of $w$ that belong to
$B$. For example if $B=\{ SEX, EXE,XES \}$ and
$w=SEXES$, then $Bad(w)$ consists of  the factors $SEX$
(occupying the first three letters), $EXE$, (occupying
letters 2,3,4), and $XES$ (occupying the last three letters).
 
So the desired generating function is also the weight-enumerator
of the much larger set consisting of pairs $(w,S)$, where
$S \subset Bad(w)$, and {\it now} the weight is defined by
$weight(w,S)=(-1)^{|S|}s^{length(w)}$. Surprisingly, it is
much easier to (weight-)count. We may think of them as
`marked words', where $S$ denotes the subset 
consisting of those words that the censor, or teacher, was able to detect.
 
First, we need a convenient data-structure for these weird objects.
Any word $w$, of length $n$, $w=w_1  \dots w_n$, has
${{n+1} \choose {2}}$ factors $w_i, \dots , w_j$, 
which we will denote by $[i,j]$.
$1 \leq i \leq j \leq n$. Hence any marked word may
be 
represented by $(w; [i_1,j_1], [i_2,j_2], \dots , [i_l, j_l])$,
where $w_{i_r} w_{i_r+1} \dots w_{j_r-1} w_{j_r} \in B$,
for $r=1, \dots , l$, and we make it canonical by
ordering the $j_r$, i.e. we arrange the marked factors
such that $j_1<j_2  < \dots < j_l$. Since no bad
word is a proper factor of another bad word, we can
assume that all the $i_r$'s are distinct, and that
there is no nesting.
 
For example if $B=\{SEX, EXE, XES \}$, and $w=SEXES$, then
$w$ gives rise to the following $2^3$ marked words:
$(SEXES; )$, $\,(SEXES;[1,3])$, $\,(SEXES;[2,4])$, 
$\,(SEXES;[3,5]),$
$\,(SEXES;[1,3],[2,4])$, 
$\,(SEXES;[1,3],[3,5]),$
$\,(SEXES;[2,4],[3,5])$, $(SEXES;[1,3],[2,4],[3,5])$.
 
For human consumption, it is easier to portray a marked
word by a 2-dimensional structure. The top line is the
word itself, and then, we list each of the factors that
are marked on a separate line, from right to left. For example the
marked word $(SEXES; )$, is simply
$$
SEXES \quad ,
$$
the marked word $(SEXES;[2,4])$, is portrayed as
$$
\matrix{ S & E & X & E & S\cr
         \space & E & X & E &\space \cr} \quad ,
$$
while the marked word $(SEXES;[1,3],[2,4],[3,5])$ is
written as:
$$
\matrix{ S & E & X & E & S\cr
  \space & \space & X & E & S \cr 
         \space & E & X & E &\space \cr 
         S & E & X & \space & \space \cr
} \quad .
$$
             
Given a word $w=w_1 \dots w_n$, 
we will say that two factors $[i,j]$ and $[i',j']$ with $j<j'$,
{\it overlap} if they have at least one common letter, i.e. if
$i<i' \leq j$.
 
Let $\M$ be the set of these marked words.
How to (weight-)count them? 
Given a non-empty marked word
$(w_1 \dots w_n; [i_1,j_1], [i_2,j_2], \dots , [i_l, j_l])$,
there are two possibilities
regarding the last letter.
 
Either $j_l<n$ , in which case $w_n$ is not part of any
detected bad factor, and deleting it results in another, shorter,
marked word
$(w_1 \dots w_{n-1}; [i_1,j_1], [i_2,j_2], \dots , [i_l, j_l])$.
We can always restore this last letter, and there is an obvious
bijection between marked words of length $n$, in which $j_l<n$ and
pairs (marked words of length $n-1$, letter of $V$).
 
The other possibility is that $j_l=n$, then we can't
simply delete the last letter $w_n$. 
Let $k$ be the smallest integer such that
$[i_k,j_k]$ overlaps with $[i_{k+1},j_{k+1}]$,
$[i_{k+1},j_{k+1} ]$ overlaps with $[i_{k+2},j_{k+2}]$,
$\dots$, $[i_{l-1},j_{l-1} ]$ overlaps with $[i_{l},j_{l}]$,
then removing the last $n-i_k+1$ letters from $w$ and
the last $l-k+1$ marked factors, results in a pair of marked words
$(w_1  \dots w_{i_k-1};[i_1,j_1], \dots , [i_{k-1},j_{k-1}])$
and $(w_{i_k} \dots w_n;[1,j_k-i_k+1], \dots , [i_l-i_k+1, j_l-i_k+1])$.
The first of these two marked words could be arbitrary, but
the second one has the special property that each of its letters
belongs to at least one marked factor, and that neighboring
marked factors overlap. Let's call such marked words
{\it clusters} and denote the set of clusters by $\C$.
 
For example if $V=\{E,S,X\}$ and
$B=\{SEX,ESE,XES\}$, the marked word 
$$
\matrix{ S & E & X & E & S & E & X\cr
  \space & \space & \space & \space & S & E & X \cr 
  \space & \space & \space & E & S & E & \space  \cr 
   S & E & X & \space & \space & \space & \space  \cr 
} \quad ,
$$
which in one-dimensional notation is $(SEXESEX;[1,3],[4,6],[5,7])$,
is not a cluster (since $[1,3]$ and $[4,6]$ don't overlap), while
$$
\matrix{ S & E & X & E & S & E & X\cr
  \space & \space & \space & \space & S & E & X \cr 
  \space & \space & \space & E & S & E & \space  \cr 
 \space & \space & X & E & S & \space & \space  \cr 
   S & E & X & \space & \space & \space & \space  \cr 
} \quad ,
$$
which in one-dimensional notation
is written $(SEXESEX;[1,3],[3,5],[4,6],[5,7])$, is a cluster.
 
Hence any member of $\M$ (i.e. marked word) is either empty
(weight $1$), or ends with a letter that is not part
of a cluster, or ends with a cluster. Peeling off
the maximal cluster, results in a smaller marked word
(by definition of maximality). Hence we have
the decomposition:
$$
\M = \{ empty\_word \} \cup \M V \cup \M\C \quad.
$$
Taking weights we have,
$$
weight(\M)= 1+ weight(\M)ds+weight(\M)weight(\C) \quad.
$$
Since $weight(\M)=f(s)$, solving for $f(s)$ yields
$$
f(s)={{1} \over {1-ds-weight(\C)}} \quad.
$$
 
It remains to find the weight-enumerator of $\C$, $weight(\C)$.
 
Let's examine how two bad words $u$ and $v$ can be the last two
members of a cluster. This happens when a proper suffix
(tail) of $u$ coincides with a proper prefix(head) of $v$.
 
For any word $w=w_1 \dots w_n$, let $HEAD(w)$ be the
set of all proper prefixes:
$$
HEAD(w_1 \dots w_n):=
\{ w_1 \,,\, w_1 w_2 \,,\, w_1 w_2 w_3 
\,,\, \dots \,,\, w_1 w_2 \dots w_{n-1}\, \} \quad ,
$$
and let $TAIL(w)$ be the set of all proper suffixes
$$
TAIL(w_1 \dots w_n):=
\{ w_n \,, \, w_{n-1} w_n \,, \, w_{n-2} w_{n-1} w_n
\, , \,\dots \,, \,w_2 \dots w_n \, \} \quad.
$$
 
Given two words $u$ and $v$, define the set
$OVERLAP(u,v):=TAIL(u) \cap HEAD(v)$. 
 
For example $OVERLAP(PICACA,CACACA)=\{CA,CACA\}$.
 
If $x \in HEAD(v)$, then we can write $v=xx'$, where $x'$ is the
word obtained from $v$ be chopping off its head $x$. Let's denote
$x'$ by $v/x$. For example $SEXYSEX/SEX=YSEX$.
 
Adopting the notation of [GrKP], section 8.4, let's define
$$
u:v \,\,:=\sum_{x \in OVERLAP(u,v)} weight(v/x) \quad,
$$
which is a certain polynomial in $s$. For example
$$
SEXSEX:EXSEXS=s+s^4 \quad,
$$
corresponding to the following two ways in which
$SEXSEX$ can be followed by $EXSEXS$ at the end of
a cluster:
$$
\matrix{ \space & E& X & S & E & X & S \cr
S & E& X & S & E & X & \space \cr
} \quad ,
$$
giving rise to weight $s$, since the leftover is the one-letter $S$,
and
$$
\matrix{ \space & \space & \space & \space & 
E& X & S & E & X & S \cr
S & E& X & S & E & X & \space & \space & \space & \space 
\cr
} \quad ,
$$
giving rise to the term $s^4$, since the leftover is the four-letter
string $SEXS$.
 
Now the set of clusters $\C$, can be partitioned into
$$
\C= \bigcup_{v \in B} \C[v] \quad ,
$$
where $\C[v]$ ($v \in B$), is the set
of clusters whose last (top) entry is $v$.
 
Given a cluster in $\C[v]$, it either consists of just $v$,
or else, chopping $v$ results in a smaller cluster that may end
with any bad word $u$ for which $OVERLAP(u,v)$ is non-empty.
This means that there is a word $x \in OVERLAP(u,v)$ for which
$u=x''x$ and $v=xx'$, for some non-empty words $x''$ and $x'$.
By removing $v$ from the cluster, we lose its tail, $x'=v/x$,
from the underlying word.
Conversely, given a cluster in $C[u]$ and one of the elements
$x$, of $OVERLAP(u,v)$, we can reconstitute the bigger cluster in
$\C[v]$ by adding $v$ to the end of the cluster, and appending the
word $v/x$ into the underlying word of the cluster.
 
For example, if once again, 
$V=\{E,S,X\}$, and $B=\{ SEX,ESE,XES\}$, then the cluster
$$
\matrix{ S & E & X & E & S & E & X\cr
  \space & \space & \space & \space & S & E & X \cr 
  \space & \space & \space & E & S & E & \space  \cr 
 \space & \space & X & E & S & \space & \space  \cr 
   S & E & X & \space & \space & \space & \space  \cr 
} \quad ,
$$
belongs to $C[SEX]$. Chopping the top $SEX$, results in the
smaller cluster
$$
\matrix{ S & E & X & E & S & E \cr
  \space & \space & \space & E & S & E & \space  \cr 
 \space & \space & X & E & S & \space & \space  \cr 
   S & E & X & \space & \space & \space & \space  \cr 
} \quad ,
$$
that belongs to $C[ESE]$, and so in this example $x=SE$, and
$x'=SEX/SE=X$.
 
We have just established a bijection
$$
\C[v] \leftrightarrow
 \{ (v,[1,length(v)]) \} \bigcup_{u \in B} \C[u] \times OVERLAP(u,v)
\quad ,
\eqno(Set\_Equations)
$$
where if $C \in \C[v]$ has more than one
bad word, and is mapped by the above bijection to $(C',x)$,
then $weight(C)=(-1)weight(C')weight(v/x)$.
 
Taking weights, we have
$$
weight(\C[v])= (-1) weight(v) - \sum_{u \in B} (u:v) \,\cdot \,
weight(\C[u]) \quad .
\eqno(Linear\_Equations)
$$
 
This is a system of $|B|$ linear equations in the $|B|$ unknowns
$weight(\C[v])$. Furthermore, it is usually rather sparse,
since for most pair of bad words  $u$ and $v$, $OVERLAP(u,v)$ is
empty. In fact, let's denote by $Comp(v)$ the set of bad words
$u \in B$ for which $OVERLAP(u,v)$ is non-empty, then the above
system can be rewritten:
$$
weight(\C[v])= -weight(v) - \sum_{u \in Comp(v)} 
(u:v) \, \cdot \, weight(\C[u]) \quad .
\eqno(Linear\_Equations')
$$
 
Note that in general $|B|$ is much smaller than $d^R$
(where $d$ is the number of letters in your alphabet $V$,
and $R+1$ is the maximal length of a bad word in $B$), the 
number of equations in the 
system of linear equations
required by the naive approach described at the beginning.
So the Goulden-Jackson method is much more efficient, in general.

After solving $(Linear\_Equations')$, we get $weight(\C)$, by using
$$
weight(\C)= \sum_{v \in B} weight(\C[v]) \quad ,
$$
which we plug into 
$$
f(s)={{1} \over {1-ds-weight(\C)}} \quad .
$$
 
{\bf Example}: Find the generating function of all words
in $\{ A, B, C, \dots , X, Y,Z \}$ that avoid the dirty
words $PIPI$ and $CACA$.
 
{\bf Answer:} $d=26$, and the system is
$$
\displaylines{
(i) \,\,weight(\C[PIPI])= -s^4-s^2 weight(\C[PIPI]) \cr
(ii) \,\,weight(\C[CACA])= -s^4-s^2 weight(\C[CACA])\cr}
$$
from which
$$
weight(\C[PIPI])=weight(\C[CACA])=-s^4/(1+s^2) \quad, 
$$
and hence $weight(\C)=-2s^4/(1+s^2)$, and hence
$$
f(s)={{1} \over {1-26s+2s^4/(1+s^2)}}=
{{1+s^2} \over {1-26s+s^2-26 s^3+2 s^4}} \quad.
$$
 
{\bf Another Example}: Find the generating function of all words
in $\{ A, B, C, \dots , X, Y,Z \}$ that avoid the dirty
words $PIPI$, $CACA$, $PICA$ and $CAPI$.
 
{\bf Answer:} $d=26$, and the system is
$$
\displaylines{
(i) \,\,weight(\C[PIPI])= -s^4-s^2 weight(\C[PIPI])-s^2 weight(\C[CAPI]) \cr
(ii) \,\,weight(\C[CACA])= -s^4-s^2 weight(\C[CACA])-s^2 weight(\C[PICA]) \cr
(iii)\,\, weight(\C[PICA])= -s^4-s^2 weight(\C[PIPI])-s^2 weight(\C[CAPI]) 
\cr
(iv)\,\, weight(\C[CAPI])= -s^4-s^2 weight(\C[CACA])-s^2 weight(\C[PICA]) 
\cr}
$$
from which
$$
f(s)=
{{1+2s^2} \over {1-26s+2s^2-52s^3+4 s^4}} \quad.
$$
\eject
{\bf Maple Implementation}
 
A Maple implementation of this
is contained in the package {\tt DAVID\_IAN}, downloadable
from this paper's website 
{\tt http://www.math.temple.edu/\Tilde zeilberg/gj.html}.
 
The function call is {\tt GJs(Alphabet,Set\_of\_bad\_words,s)}.
For example, to get the generating function 
$f(x)=\sum_{n=0}^{\infty} a(n) x^n$, where $a(n)$ is the
number of ways of spinning a dreidel $n$ times, without
having a run of length $4$ of any of Gimel, Heh, Nun, or Shin,
do \hfill \break
{\tt GJs($\{$G,H,N,S$\}$,$\{$[G,G,G,G],[H,H,H,H],[N,N,N,N],[S,S,S,S]$\}
$,x); }.
 
{\bf Penney-Ante}
 
The system of equations $(Linear\_Equations')$ is identical
to the one occurring in so-called Penney-Ante games, 
in which each player picks a word, and
a coin (or die), with as many faces as letters,
is tossed (or rolled) until a string matching that of one of the players
is encountered, in which case, she won. Since the special case of
two players and two letters is so beautifully described in [GrKP], section
8.4, and the general case is just as beautifully described in
Guibas and Odlyzko's paper [GuiO], we will only mention here how to
use our Maple implementation. The function call, in the package
{\tt DAVID\_IAN}, is {\tt Penney(List\_of\_letters,List\_of\_words,Probs)}.
The output is the list of probabilities of winning corresponding
to the list of words {\tt List\_of\_words}. {\tt Prob} is the way
the die is loaded, i.e. the probabilities of the respective
letters in the list {\tt List\_of\_letters}.
 
For example, to treat the original example in Walter Penney's
paper [P] (see also [GrKP], p. 394), in which Alice and Bob flip a coin
until either HHT or HTT occurs, and Alice wins in the former case
while Bob wins in the later case, do (in {\tt DAVID\_IAN}),
{\tt Penney([H,T],[[H,H,T],[H,T,T]],[1/2,1/2]);}, getting
the output: {\tt [2/3,1/3]}. If the probability of a Head
is $p$, then do \hfill \break
{\tt Penney([H,T],[[H,H,T],[H,T,T]],[p,1-p]);}, getting
$[p/(p^2-p+1),(1-p)^2/(p^2-p+1)]$.
 
In order to check the validity of {\tt Penney}, we have
also written a procedure {\tt PenneyGames} that simulates
many Penney-Ante games, and gives the scores of
each player. The function call is
{\tt PenneyGames(List\_of\_letters,List\_of\_words,Probs,K)},
where $K$ is the number of individual games.
Thus typing
{\tt PenneyGames([H,T],[[H,H,T],[H,T,T]],[1/2,1/2],300);}
should give something close to {\tt [200,100]}, but the
exact outcome changes, of course, for each new batch of $300$ games,
according to the whims of Lady Luck.
 
Be sure to try also {\tt BestLastPlay}, which will tell you the best
counter-move.
 
{\bf Keeping Track of the Number of Bad Words}
 
Almost nobody is perfect. It is extremely unlikely
that a long word would contain no bad factors. A more general
question is to find the number of words $a_m(n)$ in the
alphabet $V$ with {\it exactly} $m$ occurrences of factors
that belong to $B$. Let's define the generating function
$$
F(s,t):=\sum_{n=0}^{\infty}\sum_{m=0}^{n} a_m(n) s^n t^m \quad .
$$
$F(s,t)$ generalizes $f(s)$ since, obviously, $f(s)=F(s,0)$.
 
The above analysis goes almost
verbatim. Now we have:
$$
F(s,t)= \sum_{w \in V^*} weight(w) t^{\,[\,number \,\, of \,\,factors \,\, 
of \,\,w \,\, that \,\,belong \,\,to\,\, B\,]}
 \quad ,
$$
and then use the following deep facts:
$$
t= 1+ (t-1) \quad ,
$$
and for any finite set $A$,
$$
\prod_{a \in A} t = \prod_{a \in A} (1+(t-1))=
\sum_{S \subset A} (t-1)^{|S|} , \quad
$$
where as usual, $|S|$ denotes the cardinality of $S$.
 
We now have,
$$
f(s)= \sum_{w \in V^*} weight(w) t^{\,[\,number \,\, of \,\,factors \,\, 
of \,\,w \,\, that \,\,belong \,\,to\,\, B\,]}
$$
$$
= \, \sum_{w \in V^*} weight(w) (1+(t-1))
^{\,[\,number \,\, of \,\,factors \,\, 
of \,\,w \,\, that \,\,belong \,\,to\,\, B\,]}
$$
$$
= \, \sum_{w \in V^*} 
\sum_{S \subset Bad(w)} (t-1)^{|S|} s^{length(w)} \quad ,
$$
where $Bad(w)$ is the set of factors of $w$ that belong to
$B$.

The set of linear equations $(Linear\_Equations')$ now becomes:
$$
weight(\C[v])= (t-1)weight(v) + (t-1) \sum_{u \in Comp(v)} 
(u:v) \, \cdot \, weight(\C[u]) \quad ,
\eqno(Linear\_Equations'')
$$
and the rest stays the same.
 
{\bf Maple Implementation:} In the package {\tt DAVID\_IAN},
the function that finds $F(s,t)$ is {\tt GJst}. For example
let $a(n,m)$ be the number of ways of arranging $n$ children
in a line in such a way that exactly $m$ boys are isolated
(surrounded by girls on both sides, see [CoGuy], p. 205).
To find the generating
function
$F(s,t)=\sum_{n=0}^{\infty}\sum_{m=0}^{n} a(n,m)s^nt^m$
do 
{\tt GJst($\{$B,G$\}$,$\{$[G,B,G]$\}
$,s,t); }.
 
{\bf Keeping Track of the Individual Counts of Each Obscenity}
 
Suppose we want to know how many words of length $n$ has
$m_1$ occurrences of $b_1$, $m_2$ occurrences of $b_2$, \dots,
$m_f$ occurrences of $b_f$, where the set of bad words is
$B=\{ b_1, b_2, \dots b_f\}$, we need to keep track of the
individuality of each bad word. Introducing the variable
$t[b]$ for each bad word $b \in B$, we now require
$$
F(s \,;\, t[b_1] \, , \, \dots \,, \, t[b_f])= 
\sum_{w \in V^*} weight(w) \prod_{b \,\, is \,\, a \,\,bad \,\, factor
\,\,of \,\, w}
t[b]
 \quad ,
$$
and then use the following:
$$
t[b]= 1+ (t[b]-1) \quad ,
$$
and for any finite set $A$,
$$
\prod_{a \in A} t[a] = \prod_{a \in A} (1+(t[a]-1))=
\sum_{S \subset A} \prod_{a \in S} (t[a]-1)  \quad .
$$
We now have,
$$
F(s \, ; \, t[1] \, ,\, \dots \, , \, t[f] \, )= \
\sum_{w \in V^*} weight(w) 
\prod_{b \,\, is \,\, a \,\,bad \,\, factor}
t[b]
$$
$$
=\, \sum_{w \in V^*} 
\sum_{S \subset Bad(w)} 
\left ( \prod_{b \in S} (t[b]-1) \right ) s^{length(w)} \quad ,
$$
where $Bad(w)$ is the set of factors of $w$ that belong to
$B$.

The set of linear equations $(Linear\_Equations'')$ now becomes:
$$
weight(\C[v])= (t[v]-1) \cdot weight(v) + 
(t[v]-1) \cdot \sum_{u \in Comp(v)} 
(u:v) \, \cdot \, weight(\C[u]) \quad ,
\eqno(Linear\_Equations''')
$$
and the rest stays the same.
 
{\bf Maple Implementation:} In the package {\tt DAVID\_IAN},
the function that finds $F(s; t[b_1],\dots ,t[b_f] )$ 
is {\tt GJstDetail}. For example, to number of ways of arranging
$n$ kids in  line such that there are $a$ isolated boys and
$b$ isolated girls is the coefficient of $s^n t[G,B,G]^a t[B,G,B]^b$ in
the Maclaurin expansion of the rational function
{\tt GJstDetail($\{$B,G$\}$,$\{$[G,B,G], [B,G,B] $\}
$,s,t); }.
 
{\bf Keeping Track of the Letters as well}
 
If you want to know the above information, but also wish to know
the individual count of the letters, do exactly as above, with
the only difference that $weight(w)$ is no longer simply $s^{length(w)}$,
but rather (if $w=w_1 \dots w_n)$:
$$
weight(w):=\prod_{i=1}^n x[w_i] \quad .
$$
 
(For example $weight(ESSEX)=x[E]^2 x[S]^2 x[X]$). The function
calls are {\tt GJgf} and {\tt GJgfDetail}. We refer the reader to
the on-line documentation in the package {\tt DAVID\_IAN} for
instructions.
 
{\bf Generalizing to the Case of an Arbitrary Set of Bad Words}
 
What happens if we remove the condition, on the set of bad
words $B$, that no bad word can be a proper factor of another
bad word? As we saw above, if all we want is the generating
function for the number of $n$-letter words that {\it avoid}
(as factors) the members of $B$, then we can easily
remove all members of $B$ that have another member of $B$ as a factor,
until we get a set of banned words $B'$, that meets the above condition,
and that gives the same enumeration. So, as far as applying
{\tt GJs} in {\tt DAVID\_IAN}, i.e. finding the generating function
$f(s)$, we don't need to generalize.
 
{\it But} if we are interested in the more general $F(s,t)$, i.e.
in {\tt GJst}, then the original Cluster method fails. We will
now describe how to modify it.
 
Everything goes as before, but {\it now} the clusters look different.
Given a marked word\hfill\break
$(w_1 \dots w_n; [i_1,j_1], [i_2,j_2], \dots , [i_l, j_l])$,
we may no longer assume that $j_1<j_2 < \dots < j_l$, only that
$j_1 \leq j_2 \leq \dots \leq j_l$, and now we may have
nesting: i.e.: it is possible to have: $i_r<i_s<j_s<j_r$, for some $s<r$. 
Since the second component of a marked word
$(w_1 \dots w_n; [i_1,j_1], [i_2,j_2], \dots , [i_l, j_l])$
is a {\it set}, we may arrange the $[i_r,j_r]$ in such a way
that $j_r \leq j_{r+1}$ for $r=1, \dots, l-1$, and if
$j_r=j_{r+1}$, then $i_r < i_{r+1}$. For example if
$B=\{ AC,CA,CACA,ICAC,TICA,TIT,TI\}$ then
the following marked word is a cluster:
$$
\matrix{ T & I & T & I & C & A & C & A\cr
  \space & \space & \space & \space & \space & \space & C & A \cr 
  \space & \space & \space & \space & C & A & C & A \cr 
  \space & \space & \space & \space & \space & A & C & \space \cr 
  \space & \space & \space & I & C & A & C & \space \cr 
  \space & \space & T & I & C & A & \space & \space \cr 
  T & I & T & \space & \space & \space & \space & \space \cr 
  T & I & \space & \space & \space & \space & \space & \space \cr 
} \quad .
$$
In one-dimensional notation it is written:
$(TITICACA;[1,2],[1,3],[3,6],[4,7],[6,7],[5,8],[7,8])$.
 
In the original case, it was easy to enumerate clusters, since
removing the rightmost (i.e. top) bad word resulted in
a smaller cluster. This is no longer true. We are hence forced
to introduce the larger set of {\it committed clusters}.

The above marked word is a member of $\C[CA]$.
Chopping the rightmost
bad factor, $CA$, is still a cluster:
$$
\matrix{ T & I & T & I & C & A & C & A\cr
  \space & \space & \space & \space & C & A & C & A \cr 
  \space & \space & \space & \space & \space & A & C & \space \cr 
  \space & \space & \space & I & C & A & C & \space \cr 
  \space & \space & T & I & C & A & \space & \space \cr 
  T & I & T & \space & \space & \space & \space & \space \cr 
  T & I & \space & \space & \space & \space & \space & \space \cr 
} \quad ,
$$
which belongs to $\C[CACA]$,
but note that the underlying word has not changed, so  the weight
stays the same, except for a factor of $(t-1)$. If we chop the
rightmost factor again, which is now $CACA$, we get the following
cluster
$$
\matrix{ T & I & T & I & C & A & C \cr
    \space & \space & \space & \space & \space & A & C  \cr 
  \space & \space & \space & I & C & A & C  \cr 
  \space & \space & T & I & C & A & \space  \cr 
  T & I & T & \space & \space & \space & \space \cr 
  T & I & \space & \space & \space & \space & \space \cr 
} \quad ,
$$
which belongs to $\C[AC]$, BUT, unlike the previous
scenario, in which ANY cluster in $\C[AC]$ could have been gotten,
now we MUST have the $3^{rd}$ letter from the end be a $C$.
 
Such a situation occurs whenever we have $u,v \in B$ such that
$v=xuy$, where both $x$ and $y$ are {\it non-empty} words in
the alphabet $V$. For each such pair, we introduce
the set $\C'[x,u]$, which is the set of clusters whose
rightmost bad word is $u$, and the underlying word ends
with $xu$. Now we have many more unknowns and many more equations,
we set them up in an analogous way. But at the end,
after solving the system,
when we compute $weight(\C)$, we only sum
$weight(\C[v])$, and ignore all the $weight(\C'[u,x])$. Note that
$weight(\C'[u,x])$ play the roles of {\it catalysts}, that enable the
chemical reaction, but at the end are discarded.
 
We leave it to the readers to fill in the details. The readers
may get a clue from examining the Maple implementation
{\tt JODO}, that does the job, and which we will now
describe.
 
{\bf JODO: The Maple implementation of the Generalized Cluster Method}
 
The main routine is {\tt GJNZst} that computes the generating function
$F(s,t)$. For example to find the number of $10$-letter words in
the alphabet, $\{P,I\}$ containing exactly 
$13$ factors that are either $PI$, or $PIPI$, take the
coefficient of $s^{10}t^{13}$ in the Taylor expansion of
{\tt GJNZst($\{$I,P$\}$,$\{$[P,I],[P,I,P,I]$\}$,s,t);}.
 
{\bf An Interesting Application of JODO to Counting Runs}
 
A {\it run} in a word, is a string of a repeated letter.
Given a set of bad words $B$, it is of interest to know
how many words are there avoiding $B$ as factors and having
a specified number of maximal runs. It is also of interest to know
the {\it average number of maximal runs}. It can be shown that
for any finite set of bad words $B$, the average number of
runs in an $n$-letter word avoiding the words of $B$ as factors
is asymptotically $C(B)n$, where $C(B)$ is a certain algebraic number
that depends, of course, on $B$. 
 
Note that a new maximal run starts whenever we have an occurrence of
any two-letter word $ab$, with $a \neq b$. So all we have to do
is append to $B$ these words, giving them the variable $t$, and
then use a variant of {\tt GJNZ} to find the generating function.
The relevant functions are {\tt Runs} and {\tt AvRuns}.
The implementation details may be found in the package.
 
{\bf Generalizing to Non-Consecutive Bad Words}
 
So far, we wanted to avoid {\it factors}, i.e. the occurrence of a bad
word occurring as {\it consecutive} letters. Suppose we want to
avoid SEX but also the possibility that SEX would appear when
the letters are separated by one place, i.e.,
in addition to $SEX$, we don't want
factors of the form $S?E?X$, or $S?EX$, or $SE?X$,
where a question-mark could stand for any character.
Hence $SHEXY$ would be censored as would $ASELX$, but
$ASHOEOOX$ would be allowed.
In other words, we
want to include as our set of {\it bad words}, words including
a {\it blank}, where, for example, $[T,BL,T]$, means
that whenever two $T's$ are separated by {\it exactly}
one letter, we count it as a bad word. The analysis goes
almost verbatim, and the details can be found by examining the
source code of the Maple package {\tt BLANKS}, that is yet
another Maple package that accompanies this paper.
 
{\bf The Maple Package BLANKS}
 
The principal routines are {\tt BLANKSst} and 
{\tt BLANKSs0}. The function calls are
{\tt BLANKSst(alphabet, BL, MISTAKES)}
and
{\tt BLANKSs0(alphabet, BL, MISTAKES)}, where
{\tt alphabet} is the set of letters, {\tt BL} is the symbol
denoting the blank, and {\tt MISTAKES} is the set of
bad words, that are lists in the alphabet $V \cup \{BL\}$.
 
For example, to find the generating function 
$$
F(s,t):=\sum_{n}\sum_{m} a(n,m)s^n t^m  \quad ,
$$
for $a(n,m)$, the number of $0$-$1$ sequences of length $n$, 
$w=w_1, \dots, w_n$, that have
exactly $m$ occurrences of either $w_i=w_{i+1}=w_{i+2}$  or 
$w_i=w_{i+2}=w_{i+4}$ or $w_i=w_{i+3}=w_{i+6}$, type,
in {\tt BLANKS},
{\tt BLANKSst($\{$0,1$\}$,B,$\{$[0,0,0], [1,1,1], [0,B,0,B,0], [1,B,1,B,1],
[0,B,B,0,B,B,0], \hfill\break [1,B,B,1,B,B,1] $\}$)}.
 
If you want $F(s,0)$, the generating function for $a(n):=$ the
number of $0$-$1$ sequences of length $n$ with none of the
above (i.e. the number of ways of $2-$coloring the integers
$\{1,2, \dots , n \}$ such that you don't have a mono-chromatic
arithmetic sequence of length $3$ and difference $\leq 3$, then 
type:
{\tt BLANKSs0($\{1$,2$\}$,B,$\{$[0,0,0], [1,1,1], [0,B,0,B,0], [1,B,1,B,1],
[0,B,B,0,B,B,0], \hfill\break [1,B,B,1,B,B,1] $\}$)}.
 
{\bf Exploiting Symmetry}
 
Often the set of bad words is invariant either under the action of
the symmetric group (in case when the alphabet is, say, 
$\{1,2, \dots , n \}$), or under the action of the
group of signed permutations, 
(when the alphabet is, 
$\{-1,1,-2,2, \dots , -n, n \}$). Then  by symmetry,
the Cluster generating functions $weight(\C[w])$ only
depend on the equivalence class of $w$, and there are
many fewer equations, and many fewer unknowns.
The two Maple packages {\tt SYMGJ} and {\tt SPGJ}
implement these two cases respectively. We refer the readers
to the on-line documentation for details.
 
{\bf Series Expansions}
 
Many times the set of equations is too big for Maple to solve
exactly. Nevertheless, using the set of equations
$(Linear\_Equations')$ or its analogs, we can iteratively
get series expansions for the Cluster generating function, and
hence for the generating function itself, to any desired
number of terms. The procedure {\tt GJseries} in {\tt DAVID\_IAN}
handles this. The package {\tt GJseries} is a more efficient
implementation of these ideas. 
 
{\bf Applications}
 
The applications to Self-Avoiding Walks (see [MS]
for a very readable introduction to this subject) is described
in [N]. The package {\tt GJSAW}, that also comes with this
paper, is a targeted implementation.
 
Another application is to the computation of
the number of {\it ternary square-free words} (e.g. [B],[Cu]), 
which are sequences
in the alphabet $\{1,2,3\}$ that do not contain a `square' i.e.
a factor of the form $uu$ where $u$ is a word of {\it any length}.
As such, the set of bad words, $B$, is infinite, and the present
theory would have to be extended. However, we can find
{\it upper bounds} and {\it exact series expansions}, by limiting
the length of $u$. In particular, taking the
set of bad words to be $uu$, where $u$ is of length $\leq 23$,
the first $48$ terms of the sequence  
$a(n):=$ number of n-letter words in the alphabet $\{1,2,3\}$
that avoid $uu$ with $length(u) \leq 22$ coincides with the
first $46$ terms of the {\it real thing}
(i.e. $a(0)$ through $a(45)$),  and using
{\tt GJsqfree} (which is a Maple package
targeted to deal with square-free words), 
we were able to extend sequence
$M2550$ of [SP], to $46$ terms:
 
{\bf M2250} 
1, 3, 6, 12, 18, 30, 42, 60, 78, 108, 144, 204, 264, 342, 456, 618, 798,
1044, 1392, 1830, 2388, 3180, 4146, 5418, 7032, 9198, 11892, 15486, 20220,
26424, 34422, 44862, 58446, 76122, 99276, 129516, 168546, 219516, 285750,
372204, 484446, 630666, 821154, 1069512, 1392270, 1812876, 2359710, 3072486.

It is well known and easy to see (e.g. [MS], p. 9) that the
obvious inequality $a(n+m) \leq a(n)a(m)$ implies that
$\mu := \lim_{n \rightarrow \infty} a(n)^{1/n}$ exists.
 
Using Zinn-Justin's method, described in [Gut], we were able to
estimate that $\mu \approx 1.302 $, and  that if, as is reasonable
to conjecture,
$a(n) \sim \mu^n n^\theta$, then $\theta \approx 0$. 
 
Hence we have ample evidence to the following:
 
{\bf Conjecture}: The number of $n$-letter square-free
ternary words is given, asymptotically by
$a(n) \sim C \mu^n$, where 
$\mu := \lim_{n \rightarrow \infty} a(n)^{1/n}$.
 
In [B] it is shown that $2^{1/24} \approx 1.03  < \mu$,
and the upper bound $\mu \leq 1.316$ is stated.
Using the series expansion for `finite-memory' (memory $23$)
square-free words, as above, we found the sharper
upper bound $\mu \leq 1.30201064$.
 
{\bf The Maple package GJsqfree}
 
The Maple package {\tt GJsqfree}, that is also
available from this paper's website, is a targeted
implementation to the case of counting square-free words.
The main procedure is {\tt Series}, that gives the first
$NUTERMS+1$ terms of the sequence enumerating the number of
words in an alphabet of {\tt DIM} letters that avoid 
factors of the form $uu$, where the length of $u$ is
$\leq MEMO$. In particular, the first $2(MEMO+1)$ terms of this
sequence coincide with those
of the sequence of square-free words. The function call is:
{\tt Series(MEMO,DIM,NUTERMS); }. For example to get
the sequence above, we entered
{\tt Series(23,3,47); }
\bigskip
 
{\bf REFERENCES}
 
[B] J. Brinkhuis, {\it Non-repetitive sequences on three
symbols}, Quart. J. Math. Oxford (2), {\bf 34} (1983),
145-149.
 
[CoGuy], J.H. Conway and R.K Guy, {\it ``The Book of Numbers''},
Copernicus, Springer, 1996, New-York.
 
[Cu] J. Currie, {\it Open problems in pattern avoidance},
Amer. Math. Monthly {\bf 100}(1993), 790-793.
 
[GoJ1] I. Goulden  and D.M. Jackson,
{\it An inversion theorem for cluster decompositions of
sequences with distinguished subsequences},
J. London Math. Soc.(2){\bf 20} (1979), 567-576.
 
[GoJ2] I. Goulden  and D.M. Jackson,
{\it "Combinatorial Enumeration"}, John  Wiley, 1983, New York.
 
[GrKP] R.L. Graham, D.E. Knuth, and O. Patashnik,
{\it `` Concrete Mathematics''}, Addison-Welsey, 1989, Reading.
Second Edition, 1995.
 
[GuiO] L.J. Guibas and A.M. Odlyzko, {\it String overlaps,
pattern matching, and non-transitive games}, J. Comb. Theory
(ser. A) {\bf 30}(1981), 183-208.
 
[Gut]  A.J. Guttmann, {\it Asymptotic analysis of power series
expansions}, preprint.
 
[MS] N. Madras, and G. Slade, {\it ``The Self Avoiding Walk''},
Birkhauser, 1993, Boston.
 
[N] J. Noonan, {\it New upper bounds for connective constants
for self-avoiding walks}, preprint.
Available from the author's website
{\tt http://www.math.temple.edu/\Tilde noonan}.
 
[P] W. Penney, {\it Problem 95: Penney-Ante}, J. of Recreational
Math. {\bf 7}(1974), 321.
 
[SP] N.J.A. Sloane and S. Plouffe, {\it ``The Encyclopedia
of Integer Sequences''}, Academic Press, 1995.

\bye